\documentclass[11pt]{article}

\usepackage{amsmath,amssymb,amscd,epsfig,graphics}

\newtheorem{theorem}{Theorem}[section]

\newtheorem{proposition}[theorem]{Proposition}

\newtheorem{lemma}[theorem]{Lemma}

\newtheorem{remark}[theorem]{Remark}

\newcommand{\prf}{\vspace{.05in}
                    \noindent {\sc Proof} \hspace{.05in}}
\newcommand{\ethrm}{\hspace*{\fill}
                      $\Box$
                      \vspace{.1in}}

\renewcommand{\P}{{\mathbb P}}
\newcommand{\C}{{\mathbb C}}

\newcommand{\Z}{{\mathbb Z}}
\newcommand{\Q}{{\mathbb Q}}
\newcommand{\R}{{\mathbb R}}

\newcommand{\bn}{{\mathcal N}}

\newcommand{\cpl}{{\mathop{\rm cpl}\nolimits}\,}
\newcommand{\flip}{{\mathop{\rm Flip}\nolimits}}

\newcommand{\cyt}{Calabi--Yau\ threefold}

\newcommand{\toricx}{{\mathbb X}}

\newcommand{\Pdelta}{\P_\Delta}

\renewcommand{\mod}{\mathop{\rm mod}\nolimits}

\newcommand{\exc}{\mathop{\rm Exc\,}\nolimits}

\newcommand{\pic}{\mathop{\rm Pic\,}\nolimits}

\renewcommand{\phi}{\varphi}

\newcommand{\Hom}{\mathop{\rm Hom}\nolimits}
\newcommand{\spa}{\mathop{\rm span}\nolimits}

\renewcommand{\Box}{\square}

\begin{document}
 
\begin{center}

{\LARGE On a conjecture of Cox and Katz} 

\vspace{0.2in}

{\large Bal\' azs Szendr\H oi\renewcommand{\thefootnote}{}\renewcommand{\thefootnote}{$\dagger$}\footnote{Research partially supported by an Eastern European Research Bursary from Trinity College, Cambridge and an ORS Award from the British Government.}}

\vspace{0.15in}

{Mathematics Institute, University of Warwick, Coventry CV4 7AL, 

United Kingdom}

{E-mail address: \tt balazs@maths.warwick.ac.uk}

\vspace{0.15in}

{\large 13 September 2001}

\vspace{0.1in}

\end{center}

\setcounter{footnote}{0}

\renewcommand{\thefootnote}{\arabic{footnote}}

{\small
\begin{center} {\sc abstract} \end{center}
{\leftskip=30pt \rightskip=30pt
This note shows that a certain toric quotient of the quintic \cyt\ in $\P^4$ 
provides a counterexample to a recent conjecture of Cox and Katz concerning
nef cones of toric hypersurfaces.\par}}
 
\section*{Introduction}

Anticanonical hypersurfaces in toric Fano varieties provide a large supply 
of Calabi--Yau varieties. Many explicit computations and constructions on 
these varieties rely on the strength of toric geometry. In particular, ever 
since the paper of Batyrev~\cite{batyrev}, toric constructions have played an 
important role in the study of mirror symmetry. 

According to string theory, one of the basic moduli spaces involved in mirror
symmetry is the so-called K\"ahler moduli space of a Calabi--Yau variety $Z$. 
This in turn is intimately related to the nef cone $\bn(Z)$ of $Z$, i.e. the 
closed cone of divisors in the Picard group over $\R$ spanned by nef classes. 
The nef cone of a variety also appears in birational geometry; faces of this 
cone give information about possible birational contractions and fibre 
space structures on the variety. The explicit computation of this cone is
therefore often of interest.
 
For toric hypersurfaces, one could hope that this cone, or at least its 
intersection with the subspace of toric divisors on $Z$ (divisors that lie in 
the image of the restriction map from the Picard group of the ambient space), 
can be computed explicitly in terms of the toric data. In a recent 
work~\cite{cox_katz}, Cox and Katz give an combinatorial description of a 
certain cone $\bn_{0}$ in the real vector space $W$ of toric divisors on $Z$. 
This cone is constructed from the ambient toric variety and some related 
varieties birational to it; for details turn to \cite[Section 6.2]{cox_katz} 
or Section~\ref{cones} of this note. Cox and Katz conjecture 
\cite[Conjecture 6.2.8]{cox_katz} that the toric nef cone of the Calabi--Yau 
variety $Z$ is exactly $\bn_{0}$.

In this note, I consider a certain toric quotient of a quintic
hypersurface in~$\P^4$ and its toric Calabi--Yau resolution $Z$. The variety 
$Z$ is a hypersurface in a toric Fano variety and so the conjecture of 
Cox and Katz applies to it. However, I obtain 

\begin{theorem}  The hypersurface $Z$ in its ambient toric variety 
provides a counterexample to 
\cite[Conjecture~6.2.8]{cox_katz} of Cox and Katz: its (toric) nef cone is
strictly larger than the cone $\bn_{0}$ predicted by the conjecture. 
\label{extra_thm}
\end{theorem} 

My attempts at formulating an alternative conjecture or computing the nef cone
of $Z$ have not been successful. 

\vspace{0.1in}
 
\noindent
{\bf Acknowledgments} \ The particular toric quotient considered in this
note appears in a construction of Aspinwall and Morrison~\cite{asp_morr}. 
I wish to thank Pelham Wilson for encouragement, comments and corrections 
and Victor Batyrev for a discussion about nef cones.

\vspace{0.1in}
 
\noindent
{\bf Notation and conventions} \ All varieties are defined over $\C$. 
If $X$ is a projective variety, its {\it nef cone} is the closed cone
$\bn(X)$ in $\pic_\R(X)$ generated by nef divisor classes, i.e. classes
$D\in\pic(X)$ satisfying $D\cdot C\geq 0$ for all effective curves $C\subset
X$. I use the language of {\it toric geometry}, in particular the ideas of 
{\it linear Gale transform} and {\it Gelfand--Kapranov--Zelevinsky cones}; 
the main references are Fulton~\cite{fulton}, Oda--Park~\cite{oda_park} and 
Cox--Katz~\cite[Chapters~3 and 6]{cox_katz}. 

\section{Some toric varieties} 

Fix $\xi$, a primitive fifth root of unity. Consider the image $D$ of the group
\[\left\{ (z_i)\mapsto (\xi^{a_i}z_i) : \sum_{i=0}^{4} a_i = 0 \ (\mod
5)\right\}\] 
in $PGL(5,\C)$ and its subgroup $H=\langle g_1, g_2 \rangle$ generated by 
\[g_1: (z_0,z_1,z_2,z_3,z_4) \mapsto (z_0,\xi z_1,\xi^2 z_2,\xi^3 z_3,\xi^4 z_4)\] and 
\[g_2: (z_0,z_1,z_2,z_3,z_4) \mapsto (z_0,\xi z_1,\xi^3 z_2,\xi z_3,z_4).\]

I will be interested in the quotient variety $\P^4/H$ and its (partial) 
resolutions. To describe these torically, let $\widetilde N\cong \Z^4$, 
$\widetilde M=\Hom(\widetilde N,\Z)$ and consider the polyhedron 
\[
\widetilde\Delta = \left\{\sum_{i=0}^4 m_i \leq 1, \, m_i\geq -1\right\}\subset\widetilde M_\R
\]
together with its dual polyhedron 
\[\widetilde\Delta^* = \spa\{(1,0,0,0), (0,1,0,0), (0,0,1,0), (0,0,0,1),
(-1,-1,-1,-1)\}\subset \widetilde N_\R.\]
The data $(\widetilde M, \widetilde\Delta)$ defines 
$\P_{\widetilde M,\widetilde\Delta}\cong\P^4$ in the contravariant description of 
toric varieties. The obvious map of lattices and polyhedra
$(\widetilde N,\widetilde\Delta^*) \rightarrow (\widetilde M,\widetilde\Delta)$
corresponds to the quotient map 
\[\P_{\widetilde M,\widetilde\Delta} \cong \P^4   \longrightarrow \P_{\widetilde
N,\widetilde\Delta^*}\cong \P^4/D.\]

\begin{proposition} In the contravariant description, 
$\P^4/H\cong\P_{M,\Delta}$, 
where $M\cong\Z^4$ and $\Delta$ is the polyhedron
\[ \Delta = \spa\{ (1,0,0,0), (-3,5,-4,-2), (0,0,1,0), (0,0,0,1), (2,-5,3,1)\}\subset M_\R.\]
The dual polyhedron $\Delta^*\subset N_\R$ of $\Delta$ is
\[\Delta^*=\spa\{(-1,-2,-1,-1), (4,1,-1,-1), (-1,-1,-1,-1), (-1,2,4,-1), (-1,0,-1,4)\},\]
where $N=\Hom(M,\Z)$. Moreover, 
\begin{enumerate}
\item there are no lattice points in the interior of $\Delta^*$ except for
the origin;
\item there are no lattice points in the interiors of three- or 
one-dimensional faces;
\item there are precisely two lattice points $P_i, Q_i$, $i=1,\ldots, 10$ 
in the interiors of each of the ten two-dimensional faces; the combinatorics
of the faces is shown below. 
\begin{figure}[ht]\begin{center}\leavevmode\hbox{\epsfxsize=3in\epsffile{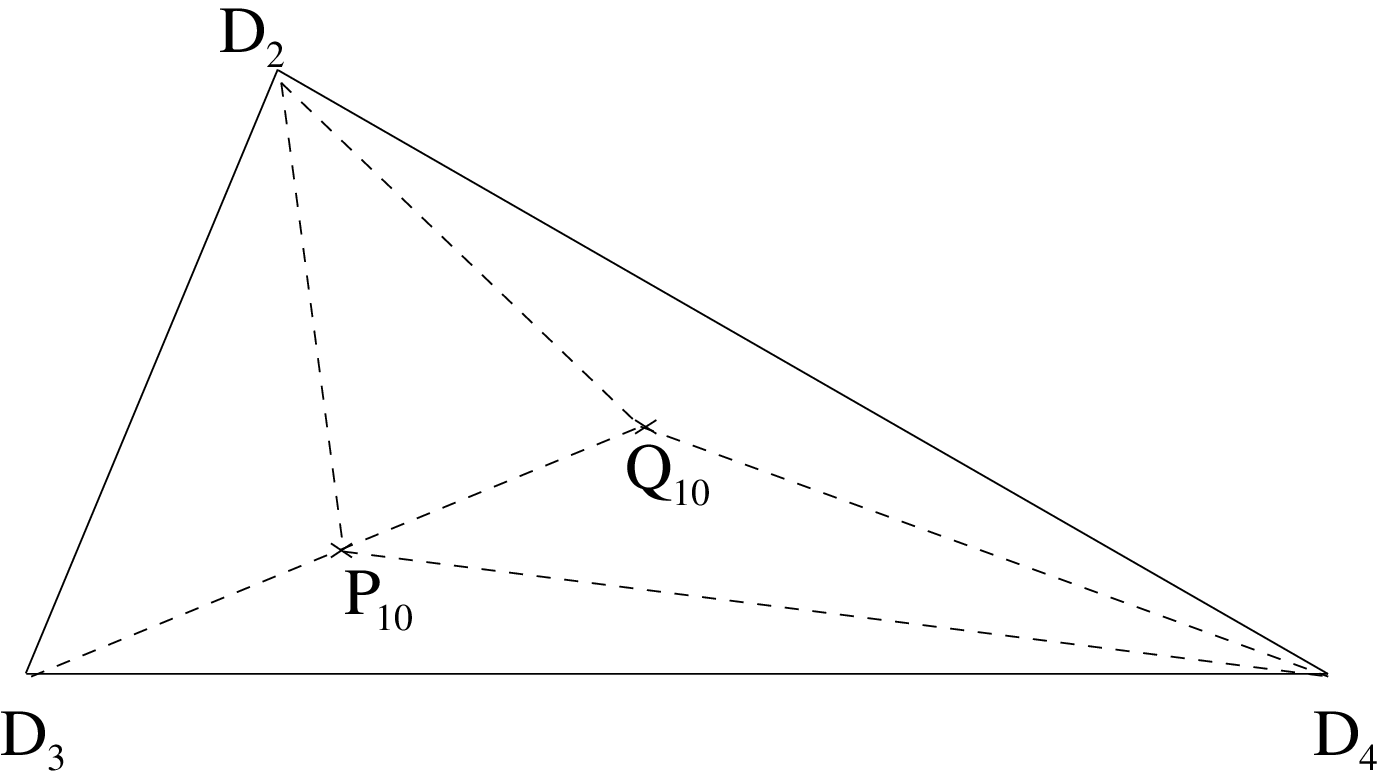}}\end{center}\label{face}\end{figure}
\end{enumerate}
\end{proposition}
\prf The inclusion $\widetilde N\hookrightarrow \widetilde M$ 
corresponds to the inclusion in 
$\widetilde M$ of the lattice of invariant monomials under the $D$ action. 
The sublattice $M$ of $\widetilde M$ is the lattice of invariant
monomials under the action of~$H$. The points $(4,-1,-1,-1)$,
$(1,-1,2,0)$, $(-1,-1,4,-1)$, $(-1,-1,-1,4)$ of $\widetilde M$ 
give a choice of basis for $M$. An easy computation gives $\Delta$, 
its dual polyhedron $\Delta^*$ and the lattice points contained in it. 
For further reference, the lattice points of $\Delta^*$ are listed in the 
Appendix. 
\ethrm

Denote $\P_{M, \Delta}$ simply by $\P_\Delta$. Let $\Sigma$ be the fan 
consisting of cones over faces of $\Delta^*$ in $N_\R$. 
This fan defines the toric variety $\toricx_{N,\Sigma}\cong \P_{\Delta}$. 

\begin{proposition} $\P_\Delta$ is a $\Q$-factorial Gorenstein variety, 
with ten curves of canonical singularities. Every permutation 
$\eta$ of the lattice points $\left\{ P_i, Q_i \right\}$ gives rise to 
a partial resolution 
\[\toricx_{\Sigma_\eta}\rightarrow\P_\Delta.\]
The varieties $\toricx_{\Sigma_\eta}$ have isolated singularities only.
\end{proposition}
\prf
All this is basic toric geometry. The curves of singularities correspond to 
the ten two-dimensional faces of $\Delta^*$. The singularities can be 
partially resolved by subdividing the fan $\Sigma$ using the lattice points 
$\left\{ P_i, Q_i \right\}$. Any permutation $\eta$ of these points  
gives a fan $\Sigma_\eta$ in the space $N_\R$ and a corresponding 
toric partial resolution $\toricx_{\Sigma_\eta}$ with isolated singularities.  
\ethrm

\section{Some hypersurfaces} 

Let $Q$ be a smooth anticanonical hypersurface in $\P^4$ invariant under the 
action of~$H$; for example the Fermat quintic will do. Let $\bar Z= Q/H$ be
the corresponding non-degenerate anticanonical hypersurface of $\P_\Delta$. 
$\bar Z$ is singular at the ten intersection points with the curves of
singularities of $\P_\Delta$ which are $\frac{1}{5}(1,1,3)$ quotient 
singularities. Every map $\toricx_{\Sigma_\eta}\rightarrow\P_\Delta$ gives
rise to a morphism $Z_\eta\rightarrow \bar Z$. The hypersurface 
$Z_\eta\subset\toricx_{\Sigma_\eta}$ is a nonsingular Calabi--Yau variety 
because every $\toricx_{\Sigma_\eta}$ is nonsingular in codimension three.

\begin{proposition} The resolutions $Z_\eta$ are all canonically isomorphic 
to a Calabi--Yau resolution $Z$ of $\bar Z$.
\end{proposition}
\prf Let $\eta_1$, $\eta_2$ be two permutations of the interior lattice 
points. There is a corresponding birational map 
\[\toricx_{\Sigma_{\eta_1}} \dashrightarrow \toricx_{\Sigma_{\eta_2}}.\]
It is easy to check that the exceptional sets of this birational map
are disjoint from the hypersurfaces $Z_{t,\eta_i}$. The statement follows. 
\ethrm

The next statement shows that in the case at hand the space of toric divisors 
is in fact the whole Picard group. 

\begin{proposition} The restriction homomorphisms 
\[\pic_\R(\toricx_{\Sigma_\eta})\longrightarrow\pic_\R(Z)\]
are all isomorphisms.   
\label{samepic}
\end{proposition}
\prf This follows from~\cite[Section 6, Theorem 2]{roan}.  
The point is that every curve of singularities of $\Pdelta$ meets the general 
hypersurface in one point, 
so the exceptional divisors in $\toricx_{\Sigma_\eta}$ 
restrict to irreducible divisors on the hypersurface $Z$. 
\ethrm

\section{Some cones}
\label{cones}

Let $W$ denote the vector space $\pic_\R(Z)$. By Proposition~\ref{samepic},
$W$ can be canonically identified with the 
{\it linear Gale transform}~\cite{oda_park} 
of the set of points $\left\{D_i, P_j, Q_r\right\}$ in~$N$.
The fans $\Sigma_\eta$ 
for different permutations $\eta$ give convex polyhedral cones, the so-called
{\it Gelfand--Kapranov--Zelevinsky cones} \ $\cpl(\Sigma_\eta)$ 
in the vector space $W$, such that there are canonical 
maps and identifications
\[\begin{array}{ccccc}\pic_\R(\toricx_{\Sigma_\eta}) & \stackrel{\sim}{\longrightarrow} & \pic_\R(Z) & = &  W \\
\cup & &&& \cup \\
\bn(\toricx_{\Sigma_\eta}) & & = && \cpl({\Sigma_\eta}) 
\end{array}\]

\begin{lemma} Under these identifications, the cones $\cpl({\Sigma_\eta})$
are all contained in the nef cone $\bn(Z)\subset W$.
\end{lemma} 
\prf The anticanonical hypersurface in $\toricx_{\Sigma_{\eta}}$ is~$Z$.
Nef divisor classes on the ambient space clearly restrict to nef classes 
on the hypersurface. 
\ethrm

Thus
\[
\bigcup_\eta \cpl({\Sigma_\eta})\subset \bn(Z).
\]
This is however not the full story. It is certainly 
possible that there are other subdivisions $\Sigma_0$ of $\Sigma$ 
satisfying the property used above; namely, that the anticanonical hypersurface
in $\toricx_{\Sigma_0}$ is isomorphic to $Z$. To treat 
these fans, I recall some definitions following~\cite[6.2]{cox_katz}. 

Suppose $\Sigma$ is a fan in $N_\R$, $\Sigma^{(1)}$ the set of its  
one-dimensional cones. A {\it linear circuit} is a linearly dependent set 
$S \subset \Sigma^{(1)}$, no subset of which is linearly dependent. 
There is a decomposition $S=S_+\cup S_-$ (depending on a choice) where 
$S_+$, respectively $S_-$ are the vectors appearing with positive, 
respectively 
negative coefficients in a linear relation. Correspondingly, there 
is a fan $\Sigma_+(S)$ given by cones spanned by $S\setminus {n_i}$ for
$n_i\in S_+$ and their subcones, and a similar fan $\Sigma_-(S)$. 

A linear circuit $S$ is said to be {\it supported by $\Sigma$}, 
if the following two conditions are satisfied: 
\begin{itemize}
\item $\Sigma_+(S)$ is a subfan of $\Sigma$. 
\item Let $\sigma$ be a top-dimensional cone of $\Sigma_+(S)$. If there 
exists a subset $S'\subset \Sigma^{(1)}$ such that $\sigma\cup S'$ generates a 
top-dimensional cone of $\Sigma$, then for all other top-dimensional 
cones $\sigma'$ of $\Sigma_+(S)$, $\sigma'\cup S'$ also 
generates a top-dimensional cone of $\Sigma$. 
\end{itemize}

\noindent 
Suppose that $S$ is a linear circuit supported
by $\Sigma$ and both  $S_-$ and $S_+$ are non-empty. 
Then there exists a new fan 
$\flip_S(\Sigma)$ from $\Sigma$ obtained by replacing the
simplices of $\Sigma$ spanned by $\sigma\cup S'$, where 
$\sigma$ is a cone of $\Sigma_-(S)$ and $S'\subset \Sigma^{(1)}$, 
by the simplices spanned
by $\sigma'\cup S'$ where $\sigma'$ is a cone of $\Sigma_+(S)$. 
Then the fan $\flip_S(\Sigma)$ is simplicial, the toric variety defined by it 
is projective, and the cones $\cpl(\Sigma)$ and $\cpl(\flip_S(\Sigma))$ in 
the linear Gale transform $W$ touch along a common face. 

Corresponding to the two fans in $N_\R$, there is a 
birational map
\[ \phi: \toricx_\Sigma \dashrightarrow \toricx_{\flip_S(\Sigma)}.
\]
It is easy to check that in case $S_-$ contains only one element, 
this map is in fact a morphism contracting a divisor. 
If however both  $S_-$ and $S_+$
contain more than one element, the birational map is a {\it generalized
flop}, a small contraction followed by a small resolution. 
If the flop $\phi$ has exceptional 
locus disjoint from the anticanonical hypersurface of
$\toricx_\Sigma$, it is called a {\it trivial flop} and in 
this case the flip attached to $S$ is referred to by~\cite{cox_katz} as a 
{\it trivial flip}. 

Return to the lattice $N$ containing the polyhedron $\Delta^*$. 
Define a fan $\Sigma_0$ in $N_\R$ to be {\it good}, if it satisfies the
following 

\vspace{0.1in}

\noindent{\bf Condition: \ }
There is a permutation $\eta$ such that the fan $\Sigma_0$ can
be obtained from the fan $\Sigma_\eta$ by a sequence of trivial flips. 

\vspace{0.1in}

The Condition implies that 
the set of one-dimensional cones of $\Sigma_0$ is precisely 
$\{D_i, P_i, Q_i\}$.
So the cones $\cpl(\Sigma_0)$ defined by the good fans $\Sigma_0$
embed canonically into~$W$. As the flips involved are trivial, the proper 
transform of~$Z$ in $\toricx_{\Sigma_0}$ is isomorphic to~$Z$. 
Setting 
\[ \bn_{0} = \bigcup_{\Sigma_0 {\ \rm good}} \cpl({\Sigma_0}),\] 
there is an inclusion
\[
\bn_{0}\subset \bn(Z).
\]
\cite[Conjecture~6.2.8]{cox_katz} expects this inclusion to be 
an equality. However, the situation is more complicated. 
The following is Theorem~\ref{extra_thm} stated in the Introduction: 

\begin{theorem} In the case discussed, the inclusion above is strict. 
The anticanonical hypersurface $Z$ in the toric variety
$\toricx_{\Sigma_{\eta}}$ provides a counterexample to 
\cite[Conjecture~6.2.8]{cox_katz} of Cox and Katz. 
\end{theorem} 

\prf Let $\Sigma_0$ be a good fan satisfying the condition that 
the cones over the tetrahedra
\[D_2P_{10}Q_{10}P_6, D_2D_{4}Q_{10}P_6, D_4P_{10}Q_{10}P_6, 
D_2P_{10}Q_{10}P_7, D_2D_{4}Q_{10}P_7, D_4P_{10}Q_{10}P_7\]
are top-di\-men\-sio\-nal cones in $\Sigma_0$ (see the figure; remember it is
a three-dimensional image of a four-dimensional setup). 

\begin{figure}[ht]\begin{center}
\leavevmode\hbox{\epsfysize=2in\epsffile{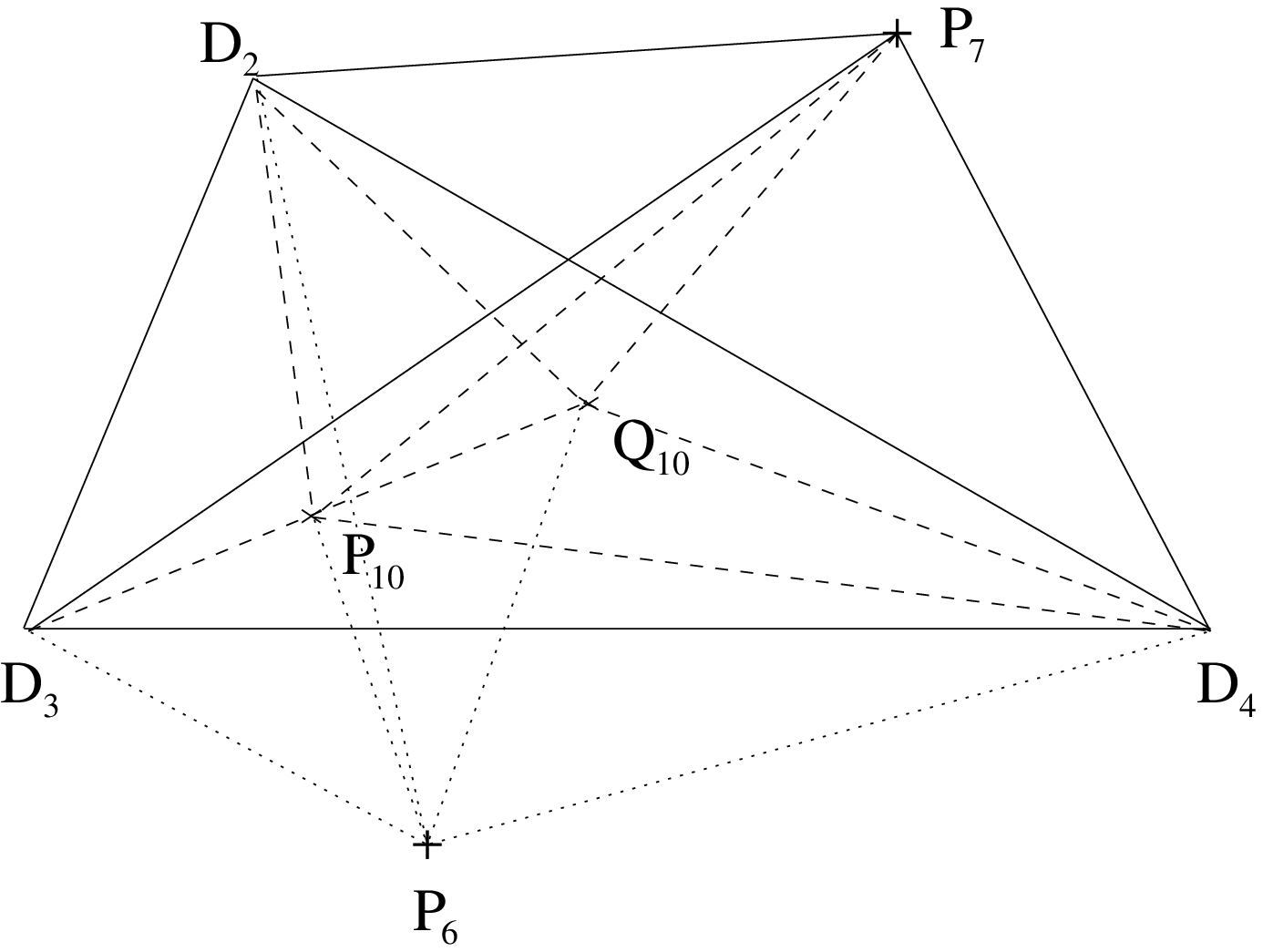}}
\end{center}\end{figure}

Under this assumption, there are two interesting circuits supported on $\Sigma_0$. The first one
is \[S=\left\{ D_2, D_4, P_{10}, Q_{10}\right\}.\] The linear relation is
\[-3{\bf q_{10}}+ {\bf d}_2 + {\bf d}_4 + {\bf p}_{10} = 0\]
in obvious notation. 
Setting $S_+= \left\{ D_2, D_4, P_{10}\right\}$, it is easy to 
check that the assumptions imply that $S$ is supported on $\Sigma_0$. 
The corresponding birational map $\phi_1$ is a contraction of the divisor $E_1$
given by the one-dimensional cone spanned by $Q_{10}$. It is easy to check that
$E_1\cong\P^1\times\P^2$ contracting to $\P^1$. $\phi_1$ restricts to the
threefold $Z$ as the contraction of an irreducible 
exceptional divisor $F_1\cong\P^2$ to a singular point. 

Now consider the circuit $S=\left\{P_6, P_7, Q_{10}\right\}$. The relation is 
\[-{\bf q}_{10}+{\bf p}_6+{\bf p}_7=0.\]
With $S_+=\left\{ P_{6}, P_7\right\}$, the circuit is supported on $\Sigma_0$. 
The corresponding contraction $\phi_2$ on $\toricx_{\Sigma_{\eta}}$ 
is again divisorial, having exceptional divisor defined by the 
one-dimensional cone spanned by $D_{10}$. So $\phi_2$ has  the same
exceptional divisor $E_2=E_1 \cong\P^1\times\P^2$ as $\phi_1$. 
The image of the exceptional divisor is in this case two-dimensional: 
the contraction $\phi_2$ restricts to $E_2$ as the projection to $\P^2$. 

The restriction of $\phi_2$ to $Z$ contracts the set 
$\exc(\phi_2)\cap Z$ to $\P^2$. 
However, in the first part of the discussion I have shown that
$\exc(\phi_2)\cap Z=\exc(\phi_1)\cap Z= \P^2$. Under $\phi_2$, this
maps isomorphically to $\P^2$. Hence the contraction $\phi_2$ 
restricts to $Z$ as the identity. 

The contraction $\phi_2$ is divisorial, in particular not a flop; 
hence the corresponding face of $\cpl(\Sigma_0)$ is a face of the 
cone $\bn_{0}$. However, divisors in this face (and beyond) are
still ample on $Z$. This implies that the corresponding face is not
in the boundary of $\bn(Z)$. Thus the cone $\bn(Z)$ is strictly 
larger than $\bn_{0}$ as claimed. 
\ethrm 

\begin{remark}\rm It is easy to see that possible counterexamples
to \cite[Conjecture~6.2.8]{cox_katz} can only arise where the
relevant face of $\bn_{0}$ 
gives a contraction with fibre dimension one. In all other cases,
the hypersurface $Z$ contains at least one contracted curve, 
and so the face is indeed a face of the nef cone of $Z$. 

From this point of 
view, it is instructive to consider the following, much simpler 
example. Let $\P=\P^1\times\P^3$ and let $Z$ be an anticanonical \cyt. 
$\P$~has a nef cone with two faces, the faces corresponding to the
contractions to the two factors. In particular, the nef cone of~$\P$ 
is also the effective cone, the cone of effective classes. 

One of the contractions restricts to $Z$ as a K3 fibration. However, 
the (Stein factorization of) the morphism to $\P^3$ 
is not a fibration, and not even a divisorial contraction: 
it is the contraction of a finite set of rational curves. 
In particular, it is a flopping face, 
there is another marked birational model for $Z$
(which as an unmarked model is incidentally isomorphic to $Z$).  
What happens here is that the nef cones are the same, but the 
effective cone changes: the effective cone of $Z$ is strictly larger
than its nef cone. Note that the trouble came again from a contraction of 
the toric ambient space of fibre dimension one. 
\end{remark} 

\section*{Appendix: Description of the polyhedron $\Delta^*$}
\label{am_vertices}

The vertices of $\Delta^*$ in $N$:
\begin{center}
$D_0  =  (-1,-2,-1,-1),\,  D_1  =  (4,1,-1,-1),\, D_2  =  (-1,-1,-1,-1),$

\vspace{0.05in}

$D_3  =  (-1,2,4,-1),\, D_4  =  (-1,0,-1,4)$  
\end{center}

\noindent The lattice points on the two-dimensional faces:
\[
\begin{array}{llcc}
P_1=(2,0,-1,-1) & Q_1=(0,-1,-1,-1)   & \mbox{ on }  &   D_0D_1D_2 \\
P_2=(0,1,2,-1) & Q_2=(1,0,0,-1) & \mbox{ on }  &  D_0D_1D_3 \\
P_3=(0,-1,-1,0) & Q_3=(1,0,-1,1)  & \mbox{ on }  &  D_0D_1D_4 \\
P_4=(-1,-1,0,-1) & Q_4=(-1,0,1,-1)   & \mbox{ on }  &  D_0D_2D_3 \\
P_5=(-1,-1,-1,0) & Q_5=(-1,-1,-1,1)   & \mbox{ on }  &  D_0D_2D_4 \\
P_6=(-1,0,0,2) & Q_6=(-1,0,1,0)   & \mbox{ on }  &  D_0D_3D_4 \\
P_7=(0,0,0,-1) & Q_7=(1,1,1,-1)  & \mbox{ on } &    D_1D_2D_3 \\
P_8=(0,0,-1,2) & Q_8=(1,0,-1,0)  & \mbox{ on }  &  D_1D_2D_4 \\
P_9=(2,1,0,0) & Q_9=(0,1,1,1) &  \mbox{ on } &  D_1D_3D_4 \\
P_{10}=(-1,1,2,0)  & Q_{10}=(-1,0,0,1) & \mbox{ on } &  D_2D_3D_4 \\
\end{array}
\]

\end{document}